\documentstyle[12pt]{article}
\def\Hom{\mathop{\rm Hom}\nolimits}
\def\End{\mathop{\rm End}\nolimits}

\def\diag{\mathop{\rm diag}\nolimits}
\def\dist{\mathop{\rm dist}\nolimits}
\def\crosspr{\mathop{%
\lefteqn{\times}\hspace{0.62em}\rule[0.01em]{0.25pt}{0.48em}\hspace{1pt}}%
\nolimits}

\def\e{{\varepsilon}}
\def\t{{\tau}}
\def\th{{\theta}}
\def\b{{\beta}}
\def\la{{\langle}}
\def\ra{{\rangle}}
\def\v{{\vert}}

\def\arr{{\longrightarrow}}
\def\a{{\alpha}}
\def\d{{\delta}}

\def\l{{\lambda}}

\def\i{{\infty}}
\def\D{{\Delta}}

\def\q{{$\quad\bullet$}}

\def\be{\begin{equation}}
\def\ee{\end{equation}}

\newtheorem{thm}{Theorem}[section]
\newtheorem{lem}[thm]{Lemma}
\newtheorem{prop}[thm]{Proposition}
\newtheorem{dfn}[thm]{Definition}

\newcommand{\norm}[1]{\left\| #1 \right\|}
\newcommand{\ov}[1]{\overline{#1}}

\def\A{{\cal A}}

\def\B{{\cal B}}

\date{2 May 1996}
\author{V.~M.~Manuilov}
\title{Diagonalizing operators over continuous fields of $C^*$-algebras}
\begin{document}

\maketitle

 \begin{abstract}
It is well known that in the commutative case, i.e. for $\A=C(X)$ being
a commutative $C^*$-algebra, compact selfadjoint operators acting on the
Hilbert $C^*$-module $H_\A$ (= continuous families of such operators
$K(x)$, $x\in X$) cannot be diagonalized inside this module but it becomes
possible if we pass to a bigger module over a bigger $W^*$-algebra
$L^\infty(X)={\bf A}\supset\A$ which can be obtained from $\A$ by
completing (on bounded sets) it with respect to the weak topology in the
natural representation of $\A$ on the Hilbert space $L^2(X)$ where the
norm is defined by a finite exact trace (measure) on $\A$.  Unlike
the ``eigenvectors'', which have coordinates from $\bf A$, the
``eigenvalues'' are continuous, i.e.  lie in the $C^*$-algebra $\A$.

We discuss here the non-commutative analog of this well-known fact. When we
pass to non-commutative $C^*$-algebras the ``eigenvalues'' are defined not
uniquely but in some cases they can also be taken from the initial
$C^*$-algebra instead of the bigger $W^*$-algebra. We prove here that such
is the case for some continuous fields of real rank zero $C^*$-algebras
over a one-dimensional manifold and give an example of a $C^*$-algebra
$\A$ for which the ``eigenvalues'' cannot be chosen from $\A$, i.e. are
discontinuous.

The main point of the proof is connected with a problem on almost commuting
operators. We prove that for some $C^*$-algebras (including the matrix
ones) if $h\in\A$ is a selfadjoint, $u\in\A$ is a unitary and if the norm
of their commutant $[u,h]$ is small enough then one can connect $u$ with
the unity by a path $u(t)$ so that the norm of the commutant $[u(t),h]$
would be also small along this path.

 \end{abstract}

\setcounter{section}{-1}
\section{Introduction}

Let $X$ be a locally compact Hausdorff space and let $\{A(x),x\in X\}$ be
a family of unital $C^*$-algebras with exact finite traces $\tau_x$,
$\tau_x(1_x)=1$.  Denote by $\prod_{x\in X}A(x)$ the set of functions
$a=a(x)$ defined on $X$ and such that $a(x)\in A(x)$ for any $x\in X$.

\begin{dfn}
Let $\A\subset
\prod_{x\in X}A(x)$ be a subset with the following properties:
\begin{enumerate}
\item
$\A$ is a $*$-subalgebra in $\prod_{x\in X}A(x)$,
\item
for any $x\in X$ the set $\{a(x),a\in\A$ is dense in the algebra $A(x)$,
\item
for any $a\in\A$ the function $x\longmapsto\norm{a(x)}$ is continuous,
\item
let $a\in\prod_{x\in X}A(x)$; if for any $x\in X$ and for any $\e>0$
one can find such $a'\in\A$ that $\norm{a(x)-a'(x)}<\e$ in some
neighborhood of the point $x$, then one has $a\in\A$,
\item
for any $a\in\A$ the function $x\longmapsto\tau_x(a(x))$ is continuous.
\end{enumerate}
Then the triple $(A(x),X,\A)$ is called a continuous field of tracial
$C^*$-algebras.
\end{dfn}

\noindent
Notice that the first four properties of this definition give the
standard definition of a continuous field of $C^*$-algebras
\cite{dix}. It is known that $\A$ is a $C^*$-algebra.
Writing in the present paper for shortness $\norm{a(x)}$ we mean the norm
of the element $a=a(x)\in\A$, i.e. $\sup_x\norm{a(x)}$.  If $V\subset X$
is a closed subset then we denote by $\A\v_V$ the restriction of this
algebra to the subspace $V$.

\smallskip\noindent
It is obvious
that $\tau:a\longmapsto\tau_x(a(x))$ is a trace on $\A$ taking values
in the $C(X)\subset\A$.

\smallskip\noindent
This trace defines an inner product on $A(x)$,
$$
(b_1,b_2)_{\t_x}=\t_x(b_1^*b_2),\qquad \norm{b}^2_{\t_x}=(b,b)_{\t_x},
$$
and a $C(X)$-valued inner product on $\A$,
\begin{equation}\label{1}
(b_1(x),b_2(x))_\t=\t_x(b_1^*(x)b_2(x))\in C(X)
\end{equation}
with the norm
$$
\norm{b(x)}_\t=\sup_x(b(x),b(x))^{1/2}_\t=\sup_x\norm{b(x)}_{\t_x}.
$$
Let $H(x)=L^2(A(x))$ be the completion of $A(x)$ with respect to the norm
$\norm{\cdot}_{\t_x}$. Then the algebra $A(x)$ is (exactly) represented
on $H(x)$ and we can pass to the corresponding $W^*$-algebra
$B(x)=L^\i(A(x))$, $A(x)\subset B(x)\subset H(x)$.

\medskip\noindent
Let now $dx$ be a $\sigma$-finite Borel measure on $X$.
Notice that the function (\ref{1}) is continuous, so one can
put
 \begin{equation}
(b_1,b_2)_\t=\int_X(b_1(x),b_2(x))_{\t_x}\,dx,\qquad
\norm{b}_\tau^2=(b,b)_\t.
 \end{equation}
Denote the completion of the
algebra $\A$ with respect to this norm by $H=L^2(\A)$.  The algebra $\A$
is (exactly) represented on $H$  and the corresponding $W^*$-algebra we
denote by ${\bf A}=L^\i(\A)\supset\A$. One can see that $$ {\bf
A}=\int^\oplus_X B(x)\,dx.  $$ From the exactness and finiteness of the
trace defined by $\tau_x$ on $\bf A$ and taking values in $L^\i(X)$ it
follows that the algebra $\bf A$ is a finite $W^*$-algebra with the finite
exact trace ${\displaystyle\overline\t=\int_X\t_x\,dx}$.

\medskip\noindent
Further on we will deal with the case when $X$ is an interval
or a circle.
In this case if $X$ is divided by points $\{x_k\}$ into smaller intervals
$D_k=[x_k;x_{k+1}]$
and if $a_k(x)\in\A\v_{D_k}$ is a continuous field for every $k$ then
we call the set $a(x)=\{a_k(x)\}$ a {\em piecewise continuous} field on
$X$.  Such piecewise continuous fields obviously belong to the
$W^*$-algebra $\bf A$.  The distance from such field to the $C^*$-algebra
$\A$ is given by the formula
 $$
\dist(a(x),\A)=\sup_k\{\norm{a_k(x_k)-a_{k+1}(x_k)}\}.
 $$
For shortness sake we write $K_1(x)\sim K_2(x)$ for two
piecewise continuous fields $K_1(x)$ and $K_2(x)$ if there exists a
piecewise continuous unitary $u(x)$ such that $u^*(x)K_1(x)u(x)=K_2(x)$.

\medskip\noindent
The present paper is organized as follows.
In the next section we discuss the notion of diagonalizability for
operators acting on Hilbert $C^*$-modules.

\smallskip\noindent
In the section \ref{a-c} we discuss a problem on almost commuting operators
and prove for some $C^*$-algebras (including the matrix ones) that if $h\in
A$ is a selfadjoint, $u\in A$ is a unitary and the norm of their commutant
$[u,h]$ is small enough then one can connect $u$ with the unity by a path
$u(t)$ so that the norm of $[u(t),h]$ would be small along this path.
We use it in the proof of the main theorem but it is also of
independent interest in view of topological applications~\cite{e-l}.
The estimate which we give for the norm of commutant is not the
best one but it is sufficient for our purpose.

\smallskip\noindent
The section \ref{theor} contains the proof of the
diagonalizability of selfadjoint compact operators on Hilbert
$C^*$-modules over some continuous fields of tracial $C^*$-algebras of
real rank zero over one-dimensional manifolds.  Unfortunately our method
cannot be applied to the case of arbitrary dimension of the base, though
we suppose that the result is still true for them. The restrictions which
we demand for the fiber $C^*$-algebras are far from necessary and can be
weakened in different ways.
In the last section we give an example of a $C^*$-algebra which does not
allow continuous diagonalization.

\medskip\noindent
{\bf Acknowledgement.}  The present paper was prepared with the partial
support of RBRF (grant N 96-01-00182). I am grateful to M.~Frank,
A.~S.~Mishchenko and E.~V.~Troitsky for helpful discussions.

\section{Diagonalizing operators on Hilbert $C^*$-mo\-dules}
\setcounter{equation}{0}

Let $\A$ be a $C^*$-algebra and
let $H_\A$ be a right Hilbert $\A$-module of sequences
$a=(a_k)$, $a_k\in \A$, $k\in {\bf N}$ such that the series $\sum a^*_ka_k$
converges in $\A$ in norm with the standard basis $\{e_k\}$ and
let $L_n(\A)\subset H_\A$ be a submodule generated by the first $n$ elements
$e_1,\ldots,e_n$ of the basis. An inner $\A$-valued product on
module $H_\A$ is given by $\la x,y\ra=\sum x^*_ky_k$ for $x,y\in\A$.
Our standard references for the theory of Hilbert $C^*$-modules and
operators on them are the
papers~\cite{pas1},\cite{pas2},\cite{rieff},\cite{kas},%
\cite{mf},\cite{lin} and the book~\cite{lance}.

\smallskip\noindent
By $H^*_\A=\Hom_\A(H_\A;\A)$ we denote the $\A$-module dual to $H_\A$
consisting of all bounded $\A$-linear $\A$-valued maps on $H_\A$.
Remember that the module $L_n(\A)$ is autodual, i.e. $L_n^*(\A)=
L_n(\A)$.
A bounded operator $K:H_\A\arr H_\A$ is
called compact~\cite{kas}~\cite{mf}, if it possesses an adjoint operator
$K^*:H_\A\arr H_\A$ and
lies in the norm closure of the
linear span of operators of the form
$$\th_{x,y},\quad \th_{x,y}(z)=x\la
y,z\ra, \quad x,y,z\in H_\A.
$$
We call a compact
operator $K$ {\it strictly positive},\/ if the operator $\la K x,x\ra$ is
positive in $\A$ and if the spectral projection corresponding to the
zero point of spectrum of $K$ is zero. We will use the well-known fact%
{}~\cite{pas2} that in the case when $\A$ is a $W^*$-algebra the
inner product can be naturally extended to the dual module $H^*_\A$.
Let $M_\A$ denote either $L_n(\A)$ or $H_\A$.

\begin{dfn}\label{diagthm}
Let $\bf A$ be a $W^*$-algebra. We call a selfadjoint operator $K$ on the
$\bf A$-module $M_{\bf A}$ {\it diagonalizable\/} if there exist a set
$\{x_i\}$ of
elements in $M^*_{\bf A}$ and a set of operators $\l_i\in\bf A$ such that
\begin{enumerate}
\item
$\{x_i\}$ is orthonormal, i.e. $\la x_i,x_j\ra=\delta_{ij}$,
\item
the module $M_{\bf A}^*$ coincides with the $\bf A$-module $\cal M^*$ dual
to the module
$\cal M$ generated by the set $\{x_i\}$,
\item
$K x_i=x_i\l_i$,
\item
for any unitaries $u_i,\,u_{i+1}\in\bf A$ we have an operator inequality
\begin{equation}\label{ord}
u_i^*\l_iu_i\ge u^*_{i+1}\l_{i+1}u_{i+1}.
\end{equation}
\end{enumerate}
\end{dfn}

\noindent
We call the elements $x_i$ ``{\it eigenvectors\/}'' and the operators $\l_i$
``{\it eigenvalues\/}'' for the operator $K$. It should be noticed that the
``eigenvectors'' and ``eigenvalues'' are defined not uniquely.
The condition (\ref{ord}) means that the notion of diagonalization
includes the natural ordering of the ``eigenvalues''.

\smallskip\noindent
The problem of diagonalizing operators in Hilbert modules was initiated by
R.~V.~Kadison in~\cite{kad} who proved that a selfadjoint operator on
the module $L_n(\bf A)$ over a $W^*$-algebra $\bf A$ can be
diagonalized. Further this problem was studied in
different settings in~\cite{gp},\cite{mur},\cite{zh},\cite{fm} etc.
In~\cite{man1},\cite{man2} we have proved the following

\begin{thm}\label{thm1}
If $\bf A$ is a finite $W^*$-algebra with a normal exact finite trace then
a compact strictly positive operator $K$ on the module $H_{\bf A}$ can be
diagonalized on $H^*_{\bf A}$ and its ``eigenvalues'' are defined uniquely
up to the unitary equivalence.  \end{thm}

\noindent
In the paper~\cite{zametki} it was proved also that the
``eigenvalues'' depend continuously on the compact operators, namely

\begin{thm}\label{contin}
If $K_r:H_{\bf A}\arr H_{\bf A}$, $r=1,\,2$
are compact
strictly positive operators and if
$\norm{K_1-K_2}<\e$ then
\begin{enumerate}
\item
one can find a unitary $U$ acting on $H_{\bf A}^*$ such that it maps the
``eigenvectors'' of $K_2$ to the ``eigenvectors'' of $K_1$ and
$\norm{U^*K_1U-K_2}<\e$,
\item
``eigenvalues'' $\{\l_i^{(r)}\}$ of operators $K_r$ $(r=1,\,2)$ can be
chosen in such a way that $\norm{\l_i^{(1)}-\l_i^{(2)}}<\e$.
\end{enumerate}
\end{thm}

\noindent
It is well known that in the commutative case, i.e. for $\A=C(X)$ being
a commutative $C^*$-algebra, compact operators cannot be diagonalized
inside $H_\A$ but it becomes possible if we pass to a bigger module over a
bigger
$W^*$-algebra $L^\infty(X)={\bf A}\supset\A$ which can be obtained from
$\A$ by completing (on bounded sets) it with respect to the weak topology
in the natural representation of $\A$ on the Hilbert space $H=L^2(\A)$ with
the new (integral) norm defined by a finite exact trace (measure) on $\A$.
Unlike the ``eigenvectors'' the ``eigenvalues'' are continuous, i.e. lie
in the $C^*$-algebra $\A$.
It leads us to the following definition.

\smallskip\noindent
Let $\A$ be a $C^*$-algebra with a finite exact trace $\t$ on it.
Let $(a,b)_\t=\t(a^*b)$, $a,b\in\cal A$, be a non-degenerate inner product
on $\A$, $\norm{\cdot}_\t=(\cdot,\cdot)^{1/2}_\t$ be the norm defined by
the trace $\t$ on $\A$.  Completing $\A$ with respect to this norm we
obtain a Hilbert space $L^2(\A)$ and an exact representation of $\A$ on
this space. Let ${\bf A}=L^\i(\A)$ be the corresponding finite tracial
$W^*$-algebra containing $\A$ as a weakly dense subalgebra.  Let $K$ be a
compact strictly positive operator on $H_\A$. We can naturally extend $K$
to the bigger module $H^*_{\bf A}$ where it will remain compact and
strictly positive and by the theorem~\ref{thm1} it can be diagonalized in
this module.

\begin{dfn}
We call a $C^*$-algebra $\A$ admitting weak
diagonalization if the diagonal entries (after diagonalization on
$H^*_{\bf A}$) for any compact
strictly positive operator $K$ on $H_\A$ can be
taken from $\A$ instead of $\bf A$.
\end{dfn}

\noindent
{\bf Problem.} Describe the class of tracial $C^*$-algebras admitting
weak diagonalization.

\medskip\noindent
Notice that if a $C^*$-algebra has the weak diagonalization property then
every selfadjoint finite rank operator $K\in M_n\otimes\A$ also can be
diagonalized over the $W^*$-algebra $\bf A$ with ``eigenvalues'' being
from $\A$.

\medskip\noindent
Recall that real rank zero ($RR(\A)=0$) means~\cite{bp} that every
selfadjoint operator in $\A$ can be approximated by operators with finite
spectrum, i.e. having the form $\sum\a_ip_i$, where
$p_i\in\A$ are selfadjoint mutually orthogonal projections and
$\a_i\in{\bf R}$. By~\cite{bp} we have in this case also
$RR(\End_\A(L_n(\A)))=0$.
In the paper~\cite{zametki} it was shown that besides the commutative
$C^*$-algebras this class contains also real rank zero tracial
$C^*$-algebras with the following property:
 \begin{itemize}
\item[($*$)]
for any two projections $p,q\in\A$ there
exist in $\A$ equivalent (in $\A$) projections $r_p\sim r_q$,
$r_p\leq p$, $r_q\leq q$ such that
$T(r_p)=T(r_q)=\min\{T(p)(z),T(q)(z)\}$,
$z\in Z$ where $L^\i(Z)$ is the center of the $W^*$-algebra
and $T$ is the standard center-valued trace on $L^\i(\A)$.
 \end{itemize}
It means that $K_0(\A)$ is
a sublattice in $K_0(L^\i(\A))$. In the case when the algebra $L^\i(\A)$
is a (finite) factor then the property ($*$) means that the map
$K_0(\A)\longrightarrow K_0(L^\i(\A))$ is a monomorphism.  Besides finite
factors this class of algebras contains the irrational rotation
$C^*$-algebras~\cite{ce} and the Bunce-Deddens algebras~\cite{bd}.

\section{On almost commuting operators}\label{a-c}
\setcounter{equation}{0}

The following proposition concerning almost commuting operators in
some $C^*$-algebras (including matrix algebras) will be used
to diagonalize continuous fields of operators. Remember that $tsr(A)=1$
means that the invertible elements are dense in $A$.

 \begin{prop}\label{a.c.}
Let $A$ be a $C^*$-algebra with properties
 \begin{enumerate}
\item
$RR(A)=0$ and $tsr(A)=1$;
\item
for every projection $p\in A$ the unitary group of the $C^*$-algebra $pAp$
is connected.
 \end{enumerate}
Let $h\in A$ be a self-adjoint and let
$u\in A$ be a unitary such that
 \be\label{ac1}
\norm{u^*hu-h}<\d.
 \ee
Then there exists a constant $C$ depending only on $\norm{h}$ and a
path $u(t)$ connecting $u$ with $1$ such that for small enough $\d$
one has for all $t$
 $$
\norm{u^*(t)hu(t)-h}<C\sqrt[4]{\d}.
 $$

\end{prop}
{\bf Proof.}
As $RR(A)=0$ we can assume without loss of generality that the operator
$h$ is a linear combination of mutually orthogonal projections,
$h=\sum_{i=1}^n \l_ip_i\in A$ with real eigenvalues $\l_i$.
We can also assume that these eigenvalues
are ordered, i.e. $\l_1>\ldots >\l_n$. Divide the segment
$[\l_n;\l_1]$ into smaller segments of the length $\sqrt[4]{\d}$ and
denote those segments which contain at least one eigenvalue $\l_i$ by
$\D_k$.  Then the number $m$ of such segments is not bigger then
$(\l_1-\l_n)/\sqrt[4]{\d}+1$ and the following properties hold:
\begin{enumerate}
\item
if $\l_i,\l_j\in\D_k$ then $\v \l_i-\l_j\v <\sqrt[4]{\d}$;
\item
if $\l_i\in\D_{k-1}$ and $\l_j\in\D_{k+1}$ then $\v \l_i-\l_j\v
\geq\sqrt[4]{\d}$.
\end{enumerate}

Let $\mu_k$ be the middle points of the segments $\D_k$. Then if
$\mu_{k+1}-\mu_k>\sqrt[4]{\d}$ then the spectrum of $h$ has a lacuna of
the length not less than $\sqrt[4]{\d}$.

\noindent
Denote the spectral projections of $h$ corresponding to the segments
$\D_k$ by $q_k$, $q_1+\ldots+q_m=1$. Then $A$ as $A$-module can be
decomposed into a direct sum corresponding to the above projections,
$A=\oplus_{k=1}^mq_kA$ and we will represent elements of the algebra $A$
as matrices with regards to this decomposition: $h=\diag(\{h_i\})$ for
$h_i=q_ihq_i$ and $u=(u_{ij})$ for $u_{ij}=q_iuq_j$.

\medskip\noindent
Notice that if $m=1$, i.e. all eigenvalues of $h$ differ from each other
not more than by $\sqrt[4]{\d}$ then we can take any path $u(t)$
connecting $u$ with $1$ because in this case there exists a number $\mu$
such that $\norm{h-\mu} <\sqrt[4]{\d}/2$, hence
 $$
\norm{v^*hv-h}\leq\norm{v^*hv-v^*\mu
v}+\norm{\mu-h}<2\frac{\sqrt[4]{\d}}{2}=\sqrt[4]{\d}
 $$
for any unitary $v$, so further we can assume that $m>1$.

\medskip\noindent
Divide once more the spectrum of the operator $h$ into smaller (than
$\D_k$) segments $\ov{\D}_s$ of length $\d$. Let $\ov{\l}_s$ be the
middle points of the segments $\ov{\D}_s$ and let $\ov{p}_s$ be the
spectral projections of $h$ corresponding to the segments $\ov{\D}_s$.
Then put
 $$
\ov{h}=\sum_s \ov{\l}_s\ov{p}_s.
 $$
Obviously $\norm{h-\ov{h}}<\d/2$, hence in view of (\ref{ac1})
 \be\label{ac.2}
\norm{u^*\ov{h}u-\ov{u}}\leq \norm{u^*\ov{h}u-u^*hu}+
\norm{u^*hu-h}+\norm{h-\ov{h}}< 2\d.
 \ee
Let $A=\oplus_s\ov{p}_sA$ be the decomposition of $A$
corresponding to the spectral projections of $h$. It is a
sub-decomposition of $\oplus_{k=1}^mq_kA$ and the matrix $u$ can be
written as $u=(v_{kl})$,
$v_{kl}=\ov{p}_k u\ov{p}_l$ and the matrix entries $u_{ij}$
can be viewed as blocks of elements $v_{kl}$.
Denote by $N$ the number of columns of the matrix $(v_{kl})$. Then one has
 \be\label{N}
N<\frac{2\norm{h}}{\d}+1
 \ee

Now we turn to construction of the homotopy $u\sim 1$. We begin with a
non-unitary path $u(t)$ which would lie not far from the set of unitaries
$U=U(A)$.

\subsection*{First step of homotopy}
For a matrix $a=(a_{ij})$ decomposed with respect to $\oplus_{k=1}^m q_kA$
we denote by $d_0(a)$ its main diagonal
$\diag(\{a_{ii}\})$, by $d_k(a)$ denote a diagonal lying $k$ lines above
(or below if $k$ is negative) the main diagonal.
We start by proving that the matrix
$u=(u_{ij})$ is ``almost'' three-diagonal.

\noindent
Let $k_i$ and $k_j$ be the numbers of the segments
$\D_k$ which contain the eigenvalues $\ov{\l}_i$ and $\ov{\l}_j$ of
$\ov{h}$ respectively; $\ov{\l}_i\in\D_{k_i}$, $\ov{\l}_j\in\D_{k_j}$.
Define a three-diagonal matrix $d(a)$
by the following way:
 \begin{enumerate}
\item
if $k_j\geq k_i +2$ or $k_j\leq k_i -2$ then put $(d(a))_{ij}=0$,
\item
if $k_j=k_i\pm 1$ and $\v\mu_{k_j}-\mu_{k_i}\v>\sqrt[4]{\d}$ then put also
$(d(a))_{ij}=0$,
\item
otherwise put $(d(a))_{ij}=a_{ij}$.
 \end{enumerate}

 \begin{lem}\label{lemma01}
It follows from (\ref{ac1}) that for small enough $\d$ one has
 $$
\norm{u-d(u)}<4\norm{h}^{1/2}\sqrt[4]{\d}.
 $$

 \end{lem}
{\bf Proof.}
Consider the matrix $(\ov{h}u-u\ov{h})(\ov{h}u-u\ov{h})^*$.
{}From the inequality
(\ref{ac.2}) it follows that the norm of this matrix
is less then $4\d^2$,
hence the norm of any element of this matrix is also less than $4\d^2$.
So, as $\ov{\l}_j$ commute with $v_{kj}$, we obtain for every $i$
 \be\label{4}
\norm{\sum_{j=1}^N (\ov{\l}_i-\ov{\l}_j)^2 v_{ij}v^*_{ij}}<4\d^2.
 \ee
Let $\ov{\l}_i\in\D_{k_i}$, $\ov{\l}_j\in\D_{k_j}$.
As all the summands in (\ref{4}) are positive, so ignoring some
of them we will not increase the norm of the sum, hence
 $$
\norm{\sum\nolimits'_j (\ov{\l}_i-\ov{\l}_j)^2 v_{ij}v^*_{ij}}<4\d^2
 $$
where the sum $\sum'$ is taken for those $j$ for which
either $\v k_j- k_i\v\geq 2$ or $\v k_j-k_i\v=1$ and
$\v\mu_{k_j}-\mu_{k_i}\v>\sqrt[4]{\d}$, i.e. we throw away those $v_{ij}$
for which $(u-d(u))_{ij}=0$.  As in the sum $\sum'$ we have
$\v\ov{\l}_k-\ov{\l}_j\v\geq\sqrt[4]{\d}$, so
 $$
4\d^2>\norm{\sum\nolimits'_j
(\ov{\l}_i-\ov{\l}_j)^2v_{ij}v^*_{ij}}\geq\sqrt{\d}\norm{\sum\nolimits'_i
v_{ij}v^*_{ij}},
 $$
hence
 \be\label{3}
\norm{\sum\nolimits'_j v_{ij}v^*_{ij}}<4\d^{3/2}.
 \ee

 \begin{lem}\label{NN}
Let $a=(a_{ij})$, $a_{ij}\in A$ be a $N\times N$ matrix such that
for every $i$ one has $\norm{\sum_j a_{ij}a^*_{ij}}<\e^2$.
Then $\norm{a}<\e\sqrt{N}$.

 \end{lem}
{\bf Proof.}
Take $\xi=(\xi_k)$, $k=1,\ldots, N$, $\xi_i\in A$.
Then using the generalized
Cauchi-Schwartz inequality~\cite{lance} we get
 \begin{eqnarray*}
\norm{a\xi}^2&=&\norm{\sum_{ijk}\xi^*_j a^*_{ij}a_{ik}\xi_k}
\leq\sum_i\norm{\sum_{kj}\xi_j^* a^*_{ij}a_{ik}\xi_k}\\
&\leq&\sum_i\norm{\sum_k a_{ik}a^*_{ik}}\cdot\norm{\sum_k \xi_k^*\xi_k}
< N\e^2 \norm{\xi}^2,
 \end{eqnarray*}
hence we have $\norm{a\xi}<\sqrt{N}\e\norm{\xi}$. \q

\smallskip\noindent
The sums $\sum'$ correspond to the blocks lying in the
matrix $u-d(u)$.
Now in view of (\ref{N}) and (\ref{3}) from the lemma \ref{NN}
(for $\e=2\d^{3/4}$) we
get for $\d<\norm{h}$
 $$
\norm{u-d(u)}<\sqrt{N}\ 2\d^{3/4}\leq
2\norm{h}^{1/2}\d^{-1/2}\,2\d^{3/4}= 4\norm{h}^{1/2}\sqrt[4]{\d}
 $$
which ends the proof of the lemma \ref{lemma01}.  \q

\medskip\noindent
Define the path $u(t)$ by the formula
 $$
u(t)=(u-d(u))(1-t)+d(u).
 $$
Then $u(0)=u$, $u(1)=d(u)$ and for every $t\in [0,1]$ by the lemma
\ref{lemma01} we have
$\dist(u(t),U)<4\norm{h}^{1/2}\sqrt[4]{\d}$.
Estimate the commutator norm:
 \begin{eqnarray*}
\norm{u(t)\ov{h}-\ov{h}u(t)}&\leq&\norm{u\ov{h}-\ov{h}u}+\norm{(u-d(u))\ov{h}-
\ov{h}(u-d(u))}\\
&<&2\d+2\norm{h}\norm{u-d(u)}<2\d+2\norm{h}4\norm{h}^{1/2}\sqrt[4]{\d}.
 \end{eqnarray*}

\subsection*{Second step of homotopy}

We should connect the matrix $d(u)$ with a diagonal one. To do so we need
the following

 \begin{lem}\label{triang}
{\rm (making matrices almost upper triangular).} \ Let
 $$
a=\left(
\begin{array}{cc}
a_{11}&a_{12}\\
a_{21}&a_{22}
\end{array}
\right)
 $$
be a matrix in $q_{j_1}A\oplus q_{j_2}A$. Then for any $\e>0$ there exists
a unitary path $v(t)$ such that $v(0)=1$ and
 $$
v(1)\cdot a=\left(
\begin{array}{cc}
a'_{11}&a'_{12}\\
a'_{21}&a'_{22}
\end{array}
\right)
 $$
with $\norm{a'_{21}}<\e$.

 \end{lem}
{\bf Proof.} As by assumption $tsr(A)=1$, so for any $\e>0$ we can
find an invertible element $\ov{a}_{11}$ such that $\norm{\ov{a}_{11}-
a_{11}}<\e$. Put $\a=a_{21}(\ov{a}_{11})^{-1}$. Put further
 $$
v(t)=\left(
\begin{array}{cc}
(1+t^2\a^*\a)^{-1/2} & (1+t^2\a^*\a)^{-1/2}\cdot t\a^*\\
-(1+t^2\a\a^*)^{-1/2}\cdot t\a & (1+t^2\a\a^*)^{-1/2}
\end{array}
\right).
 $$
It can be easily seen that $v^*(t)v(t)=v(t)v^*(t)=1$ and $v(1)$ lies in
 the unitary group in the path component of $1$. Estimate the entry
$a'_{21}$ in the product $v(1)a$:
 \begin{eqnarray*}
\norm{a'_{21}}&=&\norm{-(1+\a\a^*)^{-1/2}\a
a_{11}+(1+\a\a^*)^{-1/2}a_{21}}\\
&=&\norm{(1+\a\a^*)^{-1/2}(\a a_{11}-\a \ov{a}_{11})}
\leq \norm{(1+\a\a^*)^{-1/2}\a}\cdot\norm{a_{11}-\ov{a}_{11}}\\
&<& \e\cdot\norm{(1+\a\a^*)^{-1/2}\a}=\e\cdot
\norm{f(\b)}^{1/2}\leq\e
 \end{eqnarray*}
where $\b=\a\a^*$ and $f(\l)=\frac{\l}{1+\l}\leq 1$.\q

\medskip\noindent
Consider the element $u_{21}$
of the matrix $u$. If the segments $\D_1$ and $\D_2$ are
separated then the element $u_{21}$ is already small enough and the
element $(d(u))_{21}$ is already zero.  Put then $v_1=1$. If these
segments are not separated then by the lemma \ref{triang} we can find in
the module $q_1A\oplus q_2A$ a unitary path $v_1(t)$ with $v_1(0)=1$
and
 $$
v_1(1)=\left( \begin{array}{cc}
v^{(1)}_{11}&v^{(1)}_{12}\\
v^{(1)}_{21}&v^{(1)}_{22}
 \end{array}
\right)
 $$
such that
 $$
\left(
\begin{array}{llll}
v^{(1)}_{11}&v^{(1)}_{12}&\vdots & \\
v^{(1)}_{21}&v^{(1)}_{22}&\vdots & \\
\ldots&\ldots&\lefteqn{\ldots}\\
&&\vdots&E
\end{array}
\right)
\left(
\begin{array}{llll}
u_{11}&u_{12}&\phantom{\vdots}&\\
u_{21}&u_{22}&u_{23}\vphantom{\vdots}&\\
&\ldots&\ldots&\ldots\\
&&\ldots&\ldots\phantom{\vdots}
\end{array}
\right)
=
\left(
\begin{array}{llll}
u^{(1)}_{11}&u^{(1)}_{12}&u^{(1)}_{13}\vphantom{\vdots}&\\
u^{(1)}_{21}&u^{(1)}_{22}&u^{(1)}_{23}\vphantom{\vdots}&\\
&\ldots&\ldots&\ldots\\
&&\ldots&\ldots\phantom{\vdots}
\end{array}
\right)
 $$
with $\norm{u^{(1)}_{21}}<\e$.
Here $E$ stands for a unit matrix and empty places stand for zeros.
Notice that then
the third and the following strings of the matrix $u$ do not
change.
So we get a path $u(t)=v_1(t)d(u)$ and as $v_1(t)$ is a unitary, so
 \be\label{est}
\dist(u(t),U)<4\norm{h}^{1/2}\sqrt[4]{\d}.
 \ee
As in the lemma
\ref{triang} the number $\e$ can be taken arbitrarily small, so if we
write zero instead of the element $u^{(1)}_{21}$, the estimate (\ref{est})
would remain valid.

\noindent
Turn now to the element $u_{32}=u^{(1)}_{32}$. If the segments
$\D_2$ and $\D_3$ are separated then this element is already small enough.
Put then $v_2=1$. In the other case we can find a unitary path $v_2(t)$
with $v_2(0)=1$ and
 $$
v_2(1)=\left(
\begin{array}{ll}
v^{(2)}_{22}&v^{(2)}_{23}\\ v^{(2)}_{32}&v^{(2)}_{33}
\end{array}
\right)
 $$
in the module $q_2A\oplus q_3A$ such that
 \begin{eqnarray*}
\left(
\begin{array}{lllll}
1&&&&\\
&v^{(2)}_{22}&v^{(2)}_{23}&\vdots&\\
&v^{(2)}_{32}&v^{(2)}_{33}&\vdots&\\
&\ldots&\ldots&\lefteqn{\ldots}&\phantom{\vdots}\\
&&&\vdots&E
\end{array}
\right)
&\cdot&
\left(
\begin{array}{lllll}
u^{(1)}_{11}&u^{(1)}_{12}&u^{(1)}_{13}&\phantom{\vdots}&\\
&u^{(1)}_{22}&u^{(1)}_{23}&\phantom{\vdots}&\\
&u^{(1)}_{32}&u^{(1)}_{33}&u^{(1)}_{34}&\phantom{\vdots}\\
&&\ldots&\ldots\phantom{\vdots}&\ldots\\
&&&\ldots&\ldots\phantom{\vdots}
\end{array}
\right)
\\
&=&
\left(
\begin{array}{lllll}
u^{(2)}_{11}&u^{(2)}_{12}&u^{(2)}_{13}&&\phantom{\vdots}\\
&u^{(2)}_{22}&u^{(2)}_{23}&u^{(2)}_{24}&\phantom{\vdots}\\
&u^{(2)}_{32}&u^{(2)}_{33}&u^{(2)}_{34}&\phantom{\vdots}\\
&&\ldots&\ldots\phantom{\vdots}&\ldots\\
&&&\ldots&\ldots\phantom{\vdots}
\end{array}
\right)
 \end{eqnarray*}
with $\norm{u^{(2)}_{32}}<\e$. Again we can write zero instead of
$u^{(2)}_{32}$.
Repeating this procedure we obtain a unitary path
 $$
v(t)=v_{m-1}(t)\cdot\ldots\cdot
v_2(t)\cdot v_1(t)
 $$
connecting the unity with the unitary
 $$
v=v(m-1)=v_{m-1}(m-1)\cdot\ldots\cdot v_2(2)\cdot v_1(1).
 $$
Denote $v(t)\cdot d(u)$ by $\ov{u}(t)$ and $v\cdot d(u)$ by $\ov{u}$.
Notice that
 \be\label{dis}
\dist(\ov{u}(t),U)<4\norm{h}^{1/2}\sqrt[4]{\d}
 \ee
for all $t$. Hence we have the estimate
$\norm{\ov{u}_{ij}(t)}<4\norm{h}^{1/2}\sqrt[4]{\d}$ for every $i,j$.
Notice also that for every $t$ the path $\ov{u}(t)$ lies in the
four-diagonal matrices of the form
 $$
\ov{u}(t)=\left(
\begin{array}{llllll}
\ov{u}_{11}&\ov{u}_{12}&\ov{u}_{13}&&&\\
\ov{u}_{21}&\ov{u}_{22}&\ov{u}_{23}&\ov{u}_{24}&&\\
&\ov{u}_{32}&\ov{u}_{33}&\ov{u}_{34}&\ov{u}_{35}&\\
&&\ldots&\ldots&\ldots&\ldots\\
&&&\ldots&\ldots&\ldots\\
&&&&\ldots&\ldots
\end{array}
\right)
 $$
therefore we can easily estimate the commutator
norm along this path. To do so we should deal with the operator $h'=\sum
\mu_k q_k$, $\norm{h-h'}<\frac{1}{2}\sqrt[4]{\d}$. Then
 \begin{eqnarray}\label{com}
\norm{\ov{u}(t)h-h\ov{u}(t)}&<&\norm{[\ov{u}(t),h']}+\sqrt[4]{\d}
\nonumber\\
&\leq&\norm{[d_{-1}(\ov{u}(t)),h']}+\norm{[d_1(\ov{u}(t)),h']}+
\norm{[d_2(\ov{u}(t)),h']}+\sqrt[4]{\d}\nonumber\\
&\leq&\sup_{i,j,k}\norm{\ov{u}_{ij}(t)}(2\v\mu_k-\mu_{k+1}\v+
\v\mu_k-\mu_{k+2}\v)+\sqrt[4]{\d}\nonumber\\
&<&6\sqrt[4]{\d}
 \end{eqnarray}
for small enough $\d$. The resulting matrix $\ov{u}$ is three-diagonal and
upper triangular:
 $$
\ov{u}=
\left(
\begin{array}{lllll}
\ov{u}_{11}&\ov{u}_{12}&\ov{u}_{13}&&\\
&\ov{u}_{22}&\ov{u}_{23}&\ov{u}_{24}&\\
&&\ov{u}_{33}&\ov{u}_{34}&\ov{u}_{35}\\
&&&\ldots&\ldots\\
&&&&\ldots
\end{array}
\right).
 $$

\subsection*{Third step of homotopy}

\noindent
It follows from (\ref{dis}) that the matrix $\ov{u}$ is invertible and
there exists a unitary $w$ such that
 $$
\norm{\ov{u}-w}<4\norm{h}^{1/2}\sqrt[4]{\d}.
 $$
Therefore when $\d$ is small enough we have
 $$
\norm{(\ov{u})^{-1}-w^*}<6\norm{h}^{1/2}\sqrt[4]{\d}.
 $$
Therefore we get the estimate
 $$
\norm{(\ov{u})^{-1}-(\ov{u})^*}<10\norm{h}^{1/2}\sqrt[4]{\d},
 $$
hence
 $$
\norm{d_1(\ov{u})}<10\norm{h}^{1/2}\sqrt[4]{\d} \quad {\rm and} \quad
\norm{d_2(\ov{u})}<10\norm{h}^{1/2}\sqrt[4]{\d}.
 $$
Therefore
 \be\label{dis2}
\dist(d_0(\ov{u}),U)<\dist(\ov{u},U)+\norm{d_1(\ov{u})}+\norm{d_2(\ov{u})}
<24\norm{h}^{1/2}\sqrt[4]{\d}.
 \ee
Connect the matrices $\ov{u}$ and $d_0(\ov{u})$ by a linear path
$\ov{\ov{u}}(t)$. Then along this path one has
 \be\label{dis3}
\dist(\ov{\ov{u}}(t),U)<24\norm{h}^{1/2}\sqrt[4]{\d}
 \ee
and from (\ref{com}) we see that
 $$
\norm{\ov{\ov{u}}(t)h-h\ov{\ov{u}}(t)}<6\norm{h}^{1/2}\sqrt[4]{\d}.
 $$

\subsection*{Fourth step of homotopy}

\noindent
It follows from (\ref{dis2}) that for small enough $\d$ the diagonal
matrix $d_0(\ov{u})$ consists of invertible elements being close to
unitaries and it follows from (\ref{dis3}) that
for every $\ov{u}_{ii}$ one can find a unitary $w_i$ such that
$\norm{\ov{u}_{ii}-w_i}<24\norm{h}^{1/2}\sqrt[4]{\d}$. Take a linear path
connecting the matrices $\ov{u}$ and $\diag\{w_i\}$. Then this path also
lies close to the unitary group $U$. Then we connect the matrix
$\diag\{w_i\}$ with unity. The last path lies in $U$. Let
$\widetilde{u}(t)$ be the path of the fourth step of homotopy. Notice that
it lies within diagonal matrices, hence the commutator norm along this path
is small:
 $$
\norm{\widetilde{u}(t)h-h\widetilde{u}(t)}\leq
\norm{\widetilde{u}(t)h'-h'\widetilde{u}(t)}+\sqrt[4]{\d}=\sqrt[4]{\d}.
 $$
Consider now all four steps of constructed homotopy. We see that along the
whole path $u'(t)$ connecting $u$ with unity we have
 $$
\dist(u'(t),U)<24\norm{h}^{1/2}\sqrt[4]{\d}
 $$
and
 $$
\norm{u'(t)h-hu'(t)}<6\sqrt[4]{\d}.
 $$
Therefore there exists a unitary path $u(t)$ connecting $u$ with unity
such that
 $$
\norm{u(t)-u'(t)}<48\norm{h}^{1/2}\sqrt[4]{\d},
 $$
hence we get the estimate
 \begin{eqnarray*}
\norm{u(t)h-hu(t)}&<&\norm{u'(t)h-hu'(t)}+2\norm{h}\norm{u(t)-u'(t)}\\
&<&6\sqrt[4]{\d}+96\norm{h}^{3/2}\sqrt[4]{\d}=C\sqrt[4]{\d}
 \end{eqnarray*}
which proves the proposition. \q

\bigskip\noindent
Remark that the class of $C^*$-algebras satisfying conditions of the
proposition \ref{a.c.} includes the finite $W^*$-algebras and particularly
the matrix algebras. The constant $C$ is by no means the best one and can
be improved. On the other hand I believe that the power $1/4$ for $\d$ in
our estimate cannot be made bigger.
As the commutator norm does not change
if we replace $h$ by $h+\l$, $\l\in{\bf R}$, so the constant $C$ really
depends not on the norm of $h$ but on the length of its spectrum.

\medskip\noindent
Remark also that by the same way one can proof the similar result in the
case when instead of selfadjoint $h$ we deal with a unitary $v$ with a
lacuna in its spectrum. The last restriction is essential as otherwise
there exists an obstruction~\cite{e-l}. The following proposition can be
viewed as a first step in construction of approximately commuting homotopy
between approximately commuting and commuting pairs of unitaries in the
case when the obstruction of~\cite{e-l} is zero.

 \begin{prop}
Let $A$ be a $C^*$-algebra with properties
 \begin{enumerate}
\item
$RR(A)=0$ and $tsr(A)=1$;
\item
for every projection $p\in A$ the unitary group of the $C^*$-algebra $pAp$
is connected.
 \end{enumerate}
Let $u,v\in A$ be unitaries such that
$
\norm{u^*hu-h}<\d.
$
Suppose that the unit circle contains an arc of the length $d$ such that
the intersection of this arc with the spectrum of $v$ is empty.
Then there exists a constant $C$ depending only on $d$ and a
path $u(t)$ connecting $u$ with $1$ such that for
small enough $\d$ and for all $t$ one has
$
\norm{u^*(t)hu(t)-h}<C\sqrt[4]{\d}.
$

 \end{prop}

\section{Diagonalizing operators over some continuous fields of
$C^*$-algebras}\label{theor}
\setcounter{equation}{0}

\begin{thm}\label{wd}
Let $(A(x),X,\A)$ be a continuous field of tracial
$C^*$-algebras over an interval or a circle.
Let $(B(x),X,\B)$ be the corresponding continuous field of $W^*$-algebras
and let $\bf A$ be the corresponding $W^*$-algebra.
Suppose that for every $x\in X$ the $C^*$-algebra possesses the
following properties:
 \begin{enumerate}
\item
$RR(A(x))=0$, and $tsr(A(x))=1$;
\item
for every projection $p\in A(x)$ the unitary group of the $C^*$-algebra
$pA(x)p$ is connected;
\item
the trace $\t_x$ on $A(x)$ is finite and exact;
\item
the map $K_0(A(x))\longrightarrow K_0(B(x))$ induced by inclusion
$A(x)\subset B(x)=L^\i(A(x))$ is a monomorphism and $K_0(A(x))$ is a
sublattice in $K_0(B(x))$.
 \end{enumerate}
Then the $C^*$-algebra $\A\in \bf A$ possesses
the property of weak diagonalization.
 \end{thm}

\noindent
{\bf Proof} is based on the following lemmas.

 \begin{lem}\label{sep.lem}
Let $K_1(x)$, $K_2(x)$ be two
continuous fields of selfadjoint operators with finite spectrum in the
algebra $\A$,
$K_r$ being unitarily equivalent to $\diag \{\l^{(r)}_i\}$,
$r=1,2$. If $\norm{K_1(x) -K_2(x)}<\e$ then for any $x_0\in X$ there
exists a closed neighborhood $W$ of $x_0$ and continuous unitary fields
$u_i(x)$ on $W$ such that \be\label{estim}
\norm{\l^{(1)}_i(x)-u^*_i(x)\l^{(2)}_i(x)u_i(x)}<\e.
 \ee
 \end{lem}
{\bf Proof.}
By supposition we can write the operators $K_r(x)$ in the form
 $$
K_1(x)=\sum \a_m p_m(x),\quad K_2(x)=\sum \b_j q_j(x)
 $$
with ordered eigenvalues and $p_m(x)$, $q_j(x)$ being the spectral
projections.
As $K_0(A(x_0))\subset K_0(B(x_0))$
is a sublattice by supposition, so we can
find projections $r_l(x_0)\in A(x_0)$ such that $\sum_l r_l(x_0)=1\in
M_n\otimes A(x_0)$ and every projection $p_m(x_0)$, $q_j(x_0)$ and
projections from $M_n\otimes 1$ are unitarily equivalent to sums of some
$r_l(x_0)$.  After renumbering the eigenvalues of $K_r(x_0)$ and admitting
repeating eigenvalues we can write
 $$
K_1(x_0)=\sum\a_lr_l(x_0);\quad K_2(x_0)=\sum\b_lr'_l(x_0);\quad
r_l(x_0)\sim r'_l(x_0).
 $$
As in the $W^*$-algebra $B(x_0)=L^\i(A(x_0))$ we have
$\norm{K_1(x_0)-K_2(x_0)}<\e$ so (cf. \cite{s-t}) for all $n$ we have
 \be\label{ak}
\v\a_l-\b_l\v<\e.
 \ee
Divide the set of
projections $\{r_l(x_0)\}$ into $n$ groups
 $$
\{r_1(x_0),\ldots,r_{l_1}(x_0)\};\ldots;
 \{r_{l_{n-1}+1}(x_0),\ldots,r_{l_n}(x_0)\}
 $$
so that the sum of projections in each group would be unitarily equivalent
to the one-dimensional projection in $M_n\otimes 1$. Then each group gives
us the ``eigenvalues'' $\l^{(r)}_i(x_0)$:
 \begin{eqnarray*}
\l^{(1)}(x_0)&=&\a_{l_{i-1}+1}r_{l_{i-1}+1}(x_0)+
\ldots+\a_{l_i}r_{l_i}(x_0);\\
\l^{(2)}(x_0)&=&\b_{l_{i-1}+1}r'_{l_{i-1}+1}(x_0)+
\ldots+\b_{l_i}r'_{l_i}(x_0).
 \end{eqnarray*}
Take unitaries $u_i(x_0)\in A(x_0)$
such that
$r'_l(x_0)=u^*_i(x_0)r_l(x_0)u_i(x_0)$ for $l_{i-1}+1\leq l\leq l_i$.
Then the estimate (\ref{estim}) follows from (\ref{ak}).\q

 \begin{lem}\label{podkrutka}
Let $D=[x_k;x_{k+1}]$ be an interval in $X$ and
let $K_1(x)$, $K'_2(x)\in\A\v_D$ be two continuous fields of selfadjoint
operators with finite spectrum on $D$ with $\norm{K_1(x)-K'_2(x)}<\e$ and
let $K_2(x)\in\A(x_k)$ be such selfadjoint operator
with finite spectrum that there exists a unitary
$u(x_k)\in A(x_k)$ such that
 $$
\norm{K_2(x_k)-u^*(x_k)K'_2(x_k)u(x_k)}<\d\quad {\rm and}\quad
\norm{K_1(x_k)-K_2(x_k)}<\e.
 $$
Then for small enough $\d$ and $\e$ there
exists a piecewise continuous unitary field $u(x)$ such that on $D$ one
has
 \be\label{.0}
\norm{u^*(x)K'_2(x)u(x)-K_1(x)}<C\sqrt[4]{2\e+\d}
 \ee
where $C$ is a constant depending only on the norm of $K_1(x)$.

 \end{lem}
{\bf Proof.}
Let $u_t(x_k)$, $0\leq t\leq 1$, be a path connecting $u(x_k)$ with unity
in the unitary group of $A(x_k)$.
By assumption we have $\norm{K_2(x_k)-K'_2(x_k)}<2\e$, hence
 $$
\norm{K'_2(x_k)-u^*_t(x_k)K'_2(x_k)u_t(x_k)}<2\e+\d.
 $$
Then by the proposition \ref{a.c.}
there exists a constant $C$ depending only on the norm of $K_1(x_k)$
(or of $K'_2(x_k)$ as they are close) and
a path $u_t(x_k)$ connecting $u(x_k)$ with unity so that
 \be\label{....n}
\norm{u^*_t(x_k)K'_2(x_k)u_t(x_k)-K_1(x_k)}<C\sqrt[4]{2\e+\d}.
 \ee
It follows from (\ref{....n}) that there exists an interval $[x_k;
x'_{k+1}]$
such that for any $x\in [x_k;x'_{k+1}]$
we still have
 $$
\norm{u^*_t(x)K'_2(x)u_t(x)-K_1(x)}<C\sqrt[4]{2\e+\d}.
 $$
Put then
 $$
u(x)=u_t(x)\qquad {\rm where}
\quad \ t=\frac{x'_{k+1}-x}{x'_{k+1}-x_k}.
 $$
It is a unitary continuous field  and on $[x_k;x'_{k+1}]$
the estimate (\ref{.0})
holds. On the other hand when $x\geq x_{k+1}$ we have $u(x)=1$, so there
this estimate also holds.  \q

 \begin{lem}\label{lem1}
Let $K(x)$ be a selfadjoint operator on the module
$L_n(\A)$, $K(x)\in M_n\otimes\A$. Then for any $\e >0$ there exists a
piecewise continuous operator
field with discrete spectrum $K'(x)$ locally being from $M_n\otimes\A$
such that $\norm{K'(x)-K(x)}<\e$ and $K'(x)$ is
diagonalizable in $L_n\otimes\bf A$ with ``eigenvalues'' $\l_i(x)$ being
piecewise continuous and $\dist(\l_i(x),\A)<2\e$.
 \end{lem}

\noindent
{\bf Proof.}
Fix $\e>0$ and take a point $x_0\in X$.
As $RR(A(x_0))=0$ we can find an operator $K'(x_0)\in M_n\otimes A(x_0)$
with finite spectrum, $K'(x_0)=\sum_m\a_m p_m(x_0)$, $p_m(x_0)$ being its
spectral projections, such that $\norm{K'(x_0)-K(x_0)}<\frac{\e}{2}$.
There exists a neighborhood of the point $x_0$ such that the projections
$p_m(x_0)$ can be extended to projections $p_m(x)$ in this neighborhood.
Put $K'(x)=\sum_m\a_m p_m(x)$. Then there exists a smaller neighborhood of
$x_0$ where
 \be\label{...n}
\norm{K'(x)-K(x)}<\e.
 \ee
So every point of $X$ possesses a neighborhood $D_k$ and a continuous
field on $D_k$ such that (\ref{...n}) holds. Taking a finite
covering of $X$ by such neighborhoods we obtain a division of $X$ by
smaller intervals and a piecewise continuous field $K'(x)=\{K'_k(x)\}$,
$K'_k(x)\in M_n\otimes\A\v_{D_k}$,
such that
$\dist( K'(x),M_n\otimes\A)<2\e$. Diagonalize the operator field $K'(x)$
on every interval $D_k$, $K'(x)=\diag\{\l_i(x)\}$. Then the fields
$\l_i(x)=\{\l_i^{(k)}(x)\}$,
$\l^{(k)}_i(x)\in\A\v_{D_k}$ are piecewise continuous
and as $\norm{K'_{k-1}(x_k)-K'_k(x_k)}<2\e$ so
using lemmas
\ref{sep.lem} and \ref{podkrutka} we can change these $\l_i(x)$ on every
interval by unitarily equivalent ones to make
 $$
\dist(\l_i(x),\A)=\sup_k\norm{\l_i^{(k-1)}(x_k)-\l_i^{(k)}(x_k)}
<2\e.\qquad\bullet
 $$

\begin{lem}\label{lem2}
Let $K_1(x), K_2(x)\in M_n\otimes\bf A$ be
piecewise continuous fields of operators with finite spectrum,
$\norm{K_1(x)-K_2(x)}<\e$, $K_1(x)\sim\diag\{\l_i(x)\}$,
$K_2(x)\sim\diag\{\mu_i(x)\}$ with piecewise continuous fields $\l_i(x)$
and $\mu_i(x)$ such that $\dist(\mu_i(x),\A)<\d$.
Then there exist piecewise continuous fields of unitaries $u_i(x)$,
piecewise continuous fields $\mu'_i(x)$ and a
piecewise continuous operator field $K'_2(x)\sim\diag\{\mu'_i(x)\}$ such
that
 \begin{enumerate}
\item $\norm{K'_2(x)-K_2(x)}<\d$;
\item
$ \norm{u^*_i(x)\mu'_i(x)u_i(x)-\mu_i(x)}<\d $;
\item
$\norm{\l_i(x)-\mu'_i(x)}<C\sqrt[4]{2\e+\d}$;
\item
$\dist({\mu'_i(x),\A})<\d$.
\end{enumerate}
\end{lem}

\noindent
{\bf Proof.} Take a point $x_0\in X$ of continuity for $K_2(x)$. Let
$D\supset x_0$ be an interval in $X$.
By lemma \ref{sep.lem} we can diagonalize the operator
$K_2(x_0)$ so that its ``eigenvalues'' $\widetilde{\mu}_i(x_0)$ would
satisfy the estimate
 $$
\norm{\widetilde{\mu}_i(x_0)-\l(x_0)}<\e.
 $$
We can arbitrarily extend these ``eigenvalues'' to continuous fields
$\widetilde{\mu}_i(x)$ so that in some neighborhood of the point $x_0$
one still has
$$
\norm{\widetilde{\mu}_i(x)-\l_i(x)}<\e.
$$
By assumption there exist such unitaries $w_i\in A(x_0)$ that
$w_i^*\widetilde{\mu}(x_0)w_i=\mu_i(x_0)$. Take a unitary extensions
$w_i(x)$ of the elements $w_i$, $w_i(x_0)=w_i$ and choose a
neighborhood of the point $x_0$ so that the estimate
$$
\norm{w_i^*(x)\widetilde{\mu}_i(x)w_i(x)-\mu_i(x)}<\d
$$
would hold in this neighborhood. Then
we can get a division of $X$ by intervals
$D_k=[x_k;x_{k+1}]$ and piecewise continuous fields
$\widetilde{\mu}_i(x)$, $\widetilde{\mu}_i(x)\v_{D_k}=
\widetilde{\mu}_i^{(k)}(x)$,
$\widetilde{\mu}_i^{(k)}(x)\in\A\v_{D_k}$ and piecewise
continuous unitaries $w_i(x)=\{w_{i,k}(x)\}$,
$w_{i,k}(x)\in \A\v_{D_k}$ such that on $D_k$ one has
 \begin{equation}\label{s3}
\norm{\widetilde{\mu}_i^{(k)}(x)-\l_i(x)}<\e
 \end{equation}
and
\begin{equation}\label{s4}
\norm{w_{i,k}^*(x)\widetilde{\mu}_i^{(k)}(x)w_{i,k}(x)-\mu_i(x)}<\d/2.
\end{equation}
There exists also a piecewise continuous field of operators
$\widetilde{K}_2(x)$ unitarily equivalent to the operator
$\diag\{\widetilde{\mu}_i(x)\}$ and
 $$
\norm{\widetilde{K}_2(x)-K_2(x)}<\d/2.
 $$
It follows from (\ref{s4}) that at the point $x_k$ one has
\begin{equation}\label{s!!}
\norm{w_{i,k-1}^*(x_k)\widetilde{\mu}_i^{(k-1)}(x_k)w_{i,k-1}(x_k)-
w_{i,k}^*(x_k)\widetilde{\mu}_i^{(k)}(x_k)w_{i,k}(x_k)}<\d.
\end{equation}
It follows
from (\ref{s!!}) that there
exist such unitaries $v_{i,k}(x_k)\in A(x_k)$ that
 $$
\norm{v_{i,k}^*(x_k)\widetilde{\mu}_i^{(k-1)}(x_k)v_{i,k}(x_k)-
\widetilde{\mu}^{(k)}_i(x_k)}<\d.
 $$
By the lemma \ref{podkrutka} we can find extensions of the unitaries
$v_{i,k}(x_k)$ to some neighborhood of the point $x_k$ such that due to
(\ref{s3}) we have the estimate
 $$
\norm{v^*_{i,k}(x)\widetilde{\mu}_i^{(k)}(x)v_{i,k}(x)-\l_i(x)}<
C\sqrt[4]{2\e+\d}.
 $$
Put
 $$
\overline{\mu}_i^{(k)}(x)=v_{i,k}(x)\widetilde{\mu}_i^{(k)}(x)v_{i,k}^*(x).
 $$
Then we have
 \be\label{sdvig}
\norm{\overline{\mu}_i^{(k)}(x)-\widetilde{\mu}_i^{(k-1)}(x)}<\d
 \ee
and
 $$
\norm{\overline{\mu}_i^{(k)}(x)-\l_i(x)}<C\sqrt[4]{2\e+\d}.
 $$
Acting by induction and passing from $D_{k-1}$ to $D_k$ we change fields
$\widetilde{\mu}_i^{(k)}(x)$ by $\overline{\mu}_i^{(k)}(x)$ on every $D_k$
and so we obtain new piecewise continuous fields
$\mu'_i(x)=\{\overline{\mu}_i^{(k)}(x)\}$  such that there exist piecewise
continuous fields of unitaries $u_i(x)$ with
 $$
\norm{u^*_i(x)\mu'_i(x)u_i(x)-\mu_i(x)}<\d.
 $$
If we put $K'_2(x)=\diag\{u_i(x)\mu'_i(x)u^*_i(x)\}$ then we get validity
of $(i)$. Finally we conclude from (\ref{sdvig}) that
$\dist(\mu'_i(x),\A)<\d$.\q

\medskip\noindent
Let now $K(x)$ be a strictly positive compact
operator on the module $H_\A$. Take a sequence $\e_m>0$ converging to
zero.  Due to its compactness of $K(x)$ by the lemma~\ref{lem1} one can
find a sequence of piecewise continuous fields of diagonalizable finite
rank operators $K_m(x)\in M_{n_m}\otimes\bf A$ with ``eigenvalues''
$\l_{i,m}(x)\in\A$ such that
 \begin{enumerate}
\item
$\norm{K_m(x)-K(x)}<\e_m$,
\item
$K_m(x)\sim\diag\{\l_{i,m}(x)\}$ with $\dist(\l_{i,m}(x),\A)<\e_m$.
 \end{enumerate}
Now as
 $$
\norm{K_1(x)-K_2(x)}\leq\norm{K_1(x)-K(x)}+\norm{K(x)-K_2(x)}<\e_1+\e_2,
 $$
so by the lemma~\ref{lem2}
we can find a piecewise continuous operator
$K'_2(x)\sim\diag\{\l'_{i,2}(x)\}$ such that
 $$
\norm{K'_2(x)-K_2(x)}<\e_2,\quad
\norm{\l'_{i,2}(x)-\l_{i,2}(x)}<\e_2,\quad
\dist(\l'_{i,2}(x),\A)<\e_2,
 $$
and
 $$
\norm{\l'_{i,2}(x)-\l_{i,1}(x)}<C\sqrt[4]{2(\e_1+\e_2)+\e_2}
<2C\sqrt[4]{\e_1+\e_2}.
 $$
Then as
 \begin{eqnarray*}
\norm{K'_{m-1}(x)-K_m(x)}&\leq&\norm{K'_{m-1}(x)-K_{m-1}(x)}+
\norm{K_{m-1}(x)-K(x)}\\
&+&\norm{K(x)-K_m(x)}<2\e_{m-1}+\e_m,
 \end{eqnarray*}
so by induction we can
find a sequence $K'_m(x)\sim\diag\{\l'_{i,m}(x)\}$ such that we have
 \be\label{.1}
\norm{K'_m(x)-K_m(x)}<\e_m,
 \ee
 $$
\norm{\l'_{i,m}(x)-\l_{i,m}(x)}<\e_m,
 $$
 \be\label{.3}
\dist(\l'_{i,m}(x),\A)<\e_m,
 \ee
 \be\label{.4}
\norm{\l'_{i,m}(x)-\l_{i,m-1}(x)}<C\sqrt[4]{2(2\e_{m-1}+\e_m)+\e_m}
<2C\sqrt[4]{\e_{m-1}+\e_m}.
 \ee
As by (\ref{.1}) as
 $$
\norm{K'_m(x)-K(x)}\leq\norm{K'_m(x)-K_m(x)}+
\norm{K_m(x)-K(x)}<2\e_m,
 $$
so the sequence of operators $K'_m(x)$ converges (in norm) to the operator
$K(x)$ and from (\ref{.4}) we see that for every $i$ the sequences
$\l'_{i,m}(x)$ are Cauchi sequences provided the numbers $\e_m$ tend to
zero fast enough, namely if the series $\sum_m \sqrt[4]{\e_{m-1}+\e_m}$ is
convergent.  Hence there exist the limits
 $$
\overline{\l}_i(x)=\lim_{m\to\i}\l'_{i,m}(x)
 $$
and as by (\ref{.3})
$\dist(\l'_{i,m}(x),\A)$ tends to zero, so $\overline{\l}_i(x)\in\A$.
Show that these $\overline{\l}_i(x)$ are the ``eigenvalues'' for the
operator $K(x)$. By the theorem \ref{contin} it follows from the estimate
$\norm{K'_m(x)-K(x)}<2\e_m$ that there exist unitary operators $U_m\in
M_{n_m}\otimes\bf A$ which map the `` eigenvectors'' $x_i\in H^*_{\bf A}$
of the operator $K(x)$ to ``eigenvectors'' of the operators $K'_m(x)$ such
that $\norm{U_m^*K'_mU_m-K}<2\e_m$. Put $\widetilde{K}_m=U^*_mK'_mU_m\in
M_{n_m}\otimes\bf A$; $\widetilde{K}_m\to K$. Then one has
 \begin{equation}\label{c!}
\widetilde{K}_mx_i=x_i\l'_{i,m}.
 \end{equation}
Taking limit in (\ref{c!}) we obtain $Kx_i=x_i\overline{\l}_i$.\q

\section{Example: continuous field of rotation algebras}
\setcounter{equation}{0}

Finally we give an example of a $C^*$-algebra without property of weak
diagonalization.
Let $X=[a,b]$ be an interval of the real line containing an integer, say 1.
Let $S$ be a circle and let $C(X\times S)$ be the $C^*$-algebra of
continuous functions on the cylinder. Let $\a$ be the action of the group
$\bf Z$ of integers on this algebra defined by
 \begin{equation}\label{aa}
(\a(n)f)(x,t)=f(x,t+nx),
 \end{equation}
where $f(x,t)\in C(X\times S)$, $x\in X$, $t\in S$, $n\in\bf Z$.
Denote by
$A_X$ the crossed product
 $$
C(X\times S)\crosspr_\a\bf Z.
 $$
By~\cite{rief3},\cite{el-n-n} the algebra $A_X$ is a
continuous fields of rotation algebras $A_x=A(x)$ over an interval.
Notice that the continuous field $(A(x),X,A_X)$ is a
continuous field of tracial $C^*$-algebras and
for irrational $x$  the algebras $A(x)$ are of real rank zero and
topological rank one.
Moreover, for each $x\in X$ the algebra $A(x)$ has the property
of weak diagonalization. On the other hand the unitary group of $A(x)$ is
not connected, but the obstruction lies not here.
Notice that $A_1\cong C({\bf T}^2)$ is the commutative
$C^*$-algebra of continuous functions on a torus, hence unlike the
irrational case it is not of real rank zero and the map
$K_0(A(1))\longrightarrow K_0(B(1))$ is not a monomorphism.  There exists
in $M_2\otimes A_1$ a projection which gives the Bott generator for
$K^0({\bf T}^2)$.  This projection can be
extended~\cite{rief2},\cite{man89} to a continuous field of projections
$p(x)$ in a neighborhood of the point $1\in X$ so that for the standard
trace $\t_x$ on $A_X$ we have
 \be\label{tracex}
\t_x(p(x))=1+x.
 \ee
We
can diagonalize this field of projections in the direct integral of type
${\rm II}_1$ factors, $p(x)\sim\diag\{\l_1(x), \l_2(x)\}$ with
$\l_1(x)\geq\l_2(x)$. Then for $x=1$ we should have $\l_1(1)=1$,
$\l_2(1)=0$. For any $x$ the ``eigenvalues'' $\l_i(x)$, $i=1,2$, should be
projections. If these ``eigenvalues'' are continuous fields in $A_X$ then
for all $x$ $\l_1(x)=1$, $\l_2(x)=0$, hence
 $$
\t_x(\diag\{\l_1(x),\l_2(x)\})=1.
 $$
But as the trace is invariant under
unitary equivalence, so it contradicts (\ref{tracex}), therefore the
$C^*$-algebra $A_X$ has no property of weak diagonalization.

\smallskip\noindent
On the other hand as the algebras $A(x)$ satisfy all conditions of the
theorem \ref{wd} for a dense set of irrational points in $X$, so slightly
modifying our proof we can for given continuous field of operators $K(x)$
and for any $\e>0$ find a subset $X_\e\subset X$ so that the measure of
$X\setminus X_\e$ is less than $\e$ and $K(x)$ is diagonalizable on $X_\e$
with continuous ``eigenvalues''.

\vspace{1.5cm}
V.~M.~Manuilov\\
Dept. of Mech. and Math.,\\
Moscow State University,\\
Moscow, 119899, RUSSIA\\
e-mail: manuilov@mech.math.msu.su


\begin{thebibliography}{99}

{\small

\bibitem{bp}
{\sc L.~G.~Brown, G.~K.~Pedersen}: $C^*$-algebras of real rank zero.
{\it J. Funct. Anal.\/} \ {\bf 99} (1991), 131 -- 149.

\bibitem{bd}
{\sc J.~Bunce, J.~Deddens}: A family of simple $C^*$-algebras related to
weighted shift operators. {\it J. Funct. Anal.\/} \ {\bf 19} (1975),
13 -- 24.

\bibitem{ce}
{\sc M.-D.~Choi, G.~A.~Elliot}: Density of self-adjoint elements with
finite spectrum in an irrational rotation $C^*$-algebra. {\it Math.
Scand.\/} \ {\bf 67} (1990), 73 -- 86.

\bibitem{dix}
{\sc J.~Dixmier}: Les $C^*$-alg\`ebres et leurs repr\'esentations.
    Gauthier-Villars, Paris, 1964.

\bibitem{el-n-n}
{\sc G.~A.~Elliott, T.~Natsume, R.~Nest}: The Heisenberg group and
   $K$-theory. K\o benhavns Universitet Preprint Series No 25, 1992.

\bibitem{e-l}
{\sc R.~Exel, T.~A.~Loring}: Invariants of almost commuting unitaries.
{\it J. Funct. Anal.\/} \ {\bf 95} (1991), 364 -- 376.

\bibitem{fm}
{\sc M.~Frank, V.~M.~Manuilov}: Diagonalizing ``compact'' operators
on Hilbert $W^*$-modules. {\it Zeitschr. Anal. Anwendungen.\/} \ {\bf 14\/}
(1995), 33 --41.


\bibitem{gp}
{\sc K.~Grove, G.~K.~Pedersen}: Diagonalizing matrices over $C(X)$.
    {\it J. Funct. Anal.\/} \ {\bf 59\/} (1984), 64 -- 89.

\bibitem{kad}
{\sc R.~V.~Kadison}: Diagonalizing matrices. {\it Amer. J. Math.\/}
    \ {\bf 106\/} (1984), 1451 -- 1468.

\bibitem{kas}
{\sc G.~G.~Kasparov}: Hilbert $C^*$-modules: Theorems of Stinespring
    and Voiculescu. {\it J. Operator Theory\/} \ {\bf 4} (1980), 133
    -- 150.

\bibitem{lance}
{\sc E.~C.~Lance}: Hilbert $C^*$-modules -- a toolkit for operator
    algebraists. {\it Lecture Notes}, \/ University of Leeds, Leeds,
    1993.

\bibitem{lin}
{\sc H.~Lin}: Injective Hilbert $C^*$-modules. {\it Pacif. J.
    Math.\/} \ {\bf 154} (1992), 131 -- 164.


\bibitem{mf}
{\sc A.~S.~Mishchenko, A.~T.~Fomenko}: The index of elliptic
operators over $C^{*}$-algebras. {\it Izv. Akad. Nauk SSSR. Ser. Mat.\/} \
     {\bf 43\/} (1979), 831 -- 859 (in Russian).

\bibitem{man1}
{\sc V.~M.~Manuilov}: Diagonalization of compact operators in
     Hilbert modules over $W^{*}$-algebras of finite type. {\it Usp. Mat.
     Nauk\/} \ {\bf 49\/} (1994), No 2, 159 - 160 (in Russian).

\bibitem{man2}
{\sc V.~M.~Manuilov}: Diagonalization of compact operators in
     Hilbert modules over finite $W^{*}$-algebras. {\it Annals of Global
     Anal.  Geom.\/} \ {\bf 13\/} (1995), 207 -- 226.

\bibitem{zametki}
{\sc V.~M.~Manuilov}: Diagonalization of compact operators in
     Hilbert modules over $C^*$-algebras of real rank zero. {\it Matem
    Zametki\/}, to appear (in Russian).

\bibitem{man89}
{\sc V.~M.~Manuilov}: On $K_0$-group of continuous family of
    algebras $A_\th$. {\it Usp. Mat. Nauk\/} \ {\bf 44\/} (1989), No 3,
   163 -- 164 (in Russian).


\bibitem{mur}
{\sc Q.~J.~Murphy}: Diagonality in $C^{*}$-algebras. {\it Math.
     Zeitschr.\/} \ {\bf 199\/} (1990), 279 -- 284.

\bibitem{pas1}
{\sc W.~L.~Paschke}: Inner product modules over $B^{*}$-algebras.
     {\it Trans. Amer. Math. Soc.\/} \ {\bf 182\/} (1973), 443 -- 468.

\bibitem{pas2}
{\sc W.~L.~Paschke}: The double $B$-dual of an inner product module
     over a $C^{*}$-algebra. {\it Canad. J. Math.\/} \ {\bf 26} (1974),
   1272 -- 1280.

\bibitem{rieff}
{\sc M.~A.~Rieffel}: Induced representations of $C^*$-algebras. {\it
    Adv. in Math.\/} \ {\bf 13} (1974), 176 -- 257.

\bibitem{rief2}
{\sc M.~A.~Rieffel}: $C^*$-algebras associated with irrational
    rotations. {\it Pacif. J. Math.\/} \ {\bf 93} (1981), 415 -- 429.

\bibitem{rief3}
{\sc M.~A.~Rieffel}: Continuous fields of $C^*$-algebras coming from
group cocycles and actions. {\it Math. Ann.\/} \ {\bf 283} (1989),
631 -- 643.


\bibitem{s-t}
{\sc V.~S.~Sunder, K.~Thomsen}: Unitary orbits of selfadjoints in some
     $C^{*}$-algebras. {\it Houston J. Math.\/} {\bf 18} (1992), 127 -- 137.


\bibitem{tak}
{\sc M.~Takesaki}: Theory of operator algebras,1. New-York --
     Heidelberg -- Berlin: Springer Verlag, 1979.

\bibitem{zh}
{\sc S.~Zhang}: Diagonalizing projections in multiplier algebras and in
       matrices over a $C^*$-algebra.  {\it Pacif.  J.  Math.\/} \
       {\bf 145} (1990), 181 -- 200.

}
\end{thebibliography}
\end{document}